\documentclass[12pt]{article}
\usepackage{amsmath,amssymb}
\def\R{\mathbb{R}}
\def\C{\mathbb{C}}

\def\A{\mathbf{A}}
\def\B{\mathbf{B}}
\def\a{\mathbf{a}}
\def\b{\mathbf{b}}
\def\Be{\mathcal{B}}
\def\O{\mathcal{O}}
\def\N{\mathbb{N}}
\newtheorem{theorem}{\hspace*{\parindent}Theorem}
\newtheorem{lemma}{\hspace*{\parindent}Lemma}
\newtheorem{corollary}{\hspace*{\parindent}Corollary}

\newcounter{theremark}

\DeclareMathOperator*{\argmin}{argmin}
\title{An inverse factorial series for a general gamma ratio and related properties of the N{\o}rlund-Bernoulli polynomials}
\author{D.B.\:Karp$^{\rm a,b}$\footnote{Corresponding author. E-mail: D.B.\:Karp -- \emph{dimkrp@gmail.com},  and E.G.Prilepkina --  \emph{pril-elena@yandex.ru}}~and E.G.Prilepkina$^{\rm a,b}$
\\[10pt]
\small{\textit{$\phantom{1}^a$Far Eastern Federal University, Vladivostok, Russia}}\\\small{\textit{$\phantom{1}^b$Institute of Applied Mathematics, FEBRAS}}}
\date{}
\begin{document}
\maketitle

\begin{abstract}
We find an inverse factorial series expansion for the ratio of products of gamma functions whose arguments are linear functions of the variable.
We a give recurrence relation for the coefficients in terms of the N{\o}rlund-Bernoulli polynomials and determine quite precisely the half-plane of convergence.
Our results complement naturally a number of previous investigations of the gamma ratios which began in the 1930ies.  The expansion obtained in this paper plays a crucial role in the study of the behavior of the delta-neutral Fox's $H$ function in the neighborhood of it's finite singular point. We further apply a particular case
of the inverse factorial series expansion to derive a possibly new identity for the N{\o}rlund-Bernoulli polynomials.
\end{abstract}

\bigskip


\bigskip


\bigskip

\paragraph{1. Introduction.}  For a given complex sequence $b_0,b_1,b_2, \ldots$ the inverse factorial series is defined by
\begin{equation}\label{eq:Omega}
\Omega(z)=\sum\limits_{n=0}^{\infty}\frac{b_nn!}{(z)_{n+1}}=\sum\limits_{n=0}^{\infty}b_nB(z,n+1),
\end{equation}
where $(z)_{n+1}=z(z+1)\cdots(z+n)$ and $B(x,y)$ is Euler's beta function.  Certain variations in the denominator are possible: for instance, one can consider $(z+a)_{n+1}$ or  $z(z+w)\cdots(z+wn)$ or $\Gamma(z+a+n)/\Gamma(z+b)$ in place of $(z)_{n+1}$ for some (usually real) numbers $a,b$ and $w>0$.  This type of series had been used already by Stirling around 1730, but their rigorous theory was developed around the turn of 20th century by Landau \cite{Landau}, N{\o}rlund \cite{Norlund14,Norlund24} and Nielsen \cite{Nielsen}.  See also \cite[Chapter~X]{Milne-Thompson} for detailed account of this theory or \cite[section~4.7]{ParKam} and \cite[section~46]{Wasow} for a more concise introduction.  If the series in (\ref{eq:Omega}) converges for some value of $z$, then its domain of convergence has the form $\{z: \Re{z}>\lambda\}\!\setminus\!\{0,-1,\ldots\}$ for some real $\lambda$ called the abscissa of convergence. As the sequence $\{(z)_{n+1}^{-1}\}_{n=0}^{\infty}$ is an asymptotic sequence for $|z|\to\infty$ in $\{z:|\arg(z)|<\pi-\varepsilon\}$, the series (\ref{eq:Omega}) is an asymptotic series as $|z|\to\infty$ regardless of its convergence.  The idea to convert the Poincar\'{e} asymptotic series $\sum{a_k}z^{-k}$ into a convergent inverse factorial series goes back at least to 1912 paper of G.N.\:Watson \cite{Watson12}.  It has been revived more recently in \cite{DR} and \cite{Weniger2010} and combined with Borel-Laplace summation in the former reference.  Inverse factorial series play an important role in solution of difference equations \cite{Daalhuis,Norlund24}.

In this note we exploit similar ideas (resummation of Poincar\'{e} type series into convergent inverse factorial series) to derive the inverse factorial series expansion of the function
\begin{equation}\label{eq:W}
W(z)=\rho^{-z}\frac{\prod\nolimits_{k=1}^{p}\Gamma(A_kz+a_k)}{\prod\nolimits_{j=1}^{q}\Gamma(B_jz+b_j)}
\end{equation}
with explicit formulas for the coefficients and precise determination of the convergence domain.  Here $A_k$, $B_j$ are positive, while $a_k$, $b_j$ are complex numbers,
$\rho=\prod_{k=1}^{p}A_k^{A_k}\prod_{j=1}^{q}B_j^{-B_j}$.  This expansion was instrumental in the study of the behavior of Fox's $H$ function $H^{p,0}_{q,p}(t)$ (defined below) in the neighborhood of the singular point $t=\rho$ undertaken by us in \cite{KPCMFT2017}.  Hence, this note also fills a gap in
the proof of \cite[Theorem~1]{KPCMFT2017}.  The problem of expanding the function $W(z)$ and its particular cases in inverse factorial series has been considered previously by a number of authors.  Probably, the first appearance of such expansion is in Ford's book \cite{Ford}, where the inverse factorial series for $W(z)$ with $p=q=2$, $A_1=A_2=B_1=B_2=1/2$ and $a_1+a_2=b_1+b_2$  was found and proved to be asymptotic. No explicit formulas for the coefficients were given.  This was improved by  Van Engen in \cite{VanEngen}, where the author found the coefficients in Ford's expansion and removed the restriction $a_1+a_2=b_1+b_2$.  The general ratio $W(z)$ was first considered by Wright in the sequel \cite{Wright1,Wright2}. He proved that there exists a series in reciprocal gamma functions asymptotic to the function $W(z)$ under very general assumptions.  Wright only gave a formula for the first coefficient, but mentioned that further coefficients could also be computed.  Similar result was later proved by Hughes in \cite{Hughes45} under the natural additional restriction $\sum_{k=1}^{p}A_k=\sum_{j=1}^{q}B_j$ and using the standard inverse factorial series (\ref{eq:Omega}).   In his milestone work \cite{Norlund55} N{\o}rlund deduced an inverse factorial series expansion of the function (\ref{eq:W}) when $p=q$ and $A_k=B_k=1$, $k=1,\ldots,p$, and proved its convergence in the intersection of the half planes $\Re(z+a_k)>0$, $k=1,\ldots,p$.  He also gave two different methods to compute the coefficients in this expansion.  In a series of papers \cite{Riney1}-\cite{Riney3} Riney studied the function (\ref{eq:W}) for $p\leq{q+1}$ and $A_j=B_k=1$, $j=1,\ldots,p$, $k=1,\ldots,q$.  He gave an asymptotic series for this function in terms of gamma ratios, of which standard factorial series is a particular case, and presented several methods to compute the coefficients. Riney's investigations were complemented by van der Corput \cite{vanCorput}, who considered the opposite case $q\leq{p+1}$, and Wright \cite{Wright3}, who suggested further methods for calculating the coefficients.  Braaksma \cite{Braaksma} again considered the general case of (\ref{eq:W}) and proved that there exists a series in reciprocal gamma functions asymptotic to $W(z)$. He also gave an explicit formula for the principal term. This result of Braaksma (which is just  a technical tool in his deep investigation of Mellin-Barnes integrals) is, in fact, a modification of the earlier work by Wright \cite{Wright1,Wright2} mentioned above.  A survey of some of the above work is given in section~2.2 of the book \cite{ParKam} by Paris and Kaminski, where one can also find explicit proofs and several examples. Independently, Gupta and Tang \cite{GuptaTang} presented a series in gamma ratios for $W(z)$ when  $\sum_{k=1}^{p}A_k=\sum_{j=1}^{q}B_j$ and gave certain recursions for computing the coefficients. They also claim convergence but gave no real proof of this claim.  Further details about their work can be found in the introduction to our paper \cite{KPCMFT2017}.

 In the present paper we combine some ideas from \cite{GuptaTang}  with Borel-Laplace summation to give a rigorous proof of convergence and formulas for the coefficients for the inverse factorial series expansion of $W(z)$ under the assumption $\sum_{k=1}^{p}A_k=\sum_{j=1}^{q}B_j$.  This is done in section~3 of this paper. Furthermore, in section~4 we apply the inverse factorial series for a simplest particular case of (\ref{eq:Omega}) to derive a presumably new identity for the N{\o}rlund-Bernoulli polynomials. The main results presented in sections~3 and 4 are preceded by the preliminaries expounded in section~2.

\paragraph{2. Preliminaries.}  A function $f(z)$ of a complex variable $z$ is said to possess a Poincar\'{e} type asymptotic expansion in an unbounded domain $D$ if
\begin{equation}\label{eq:asimpt}
f(z)=\sum\limits_{k=0}^{n-1}\frac{\alpha_k}{z^k}+R_n(z)
\end{equation}
and $R_n(z)=\O(z^{-n})$ as $D\ni{z}\to\infty$. If this holds for all natural $n$ it is customary to write
\begin{equation*}
f(z)\sim\sum\limits_{k=0}^{\infty}\frac{\alpha_k}{z^k}~~~\text{as}~z\to\infty~\text{in}~D.
\end{equation*}
The asymptotic expansion of $f(z)$ is said to be Gevrey-1 (or Gevrey of order 1),
if there exist the numbers $M,\tau>0$ such that for all $z\in{D}$ and all positive integers $n$ the error term $R_n(z)$ in (\ref{eq:asimpt}) satisfies  \cite[Definition~5.21]{Mitschi-Sauzin}, \cite[Definition~4.130]{Cost}
\begin{equation}\label{eq:remainder}
|R_n(z)|\leq\frac{M\tau^n n!}{|z|^n}.
\end{equation}
It is convenient to introduce the following class of functions.

\smallskip

\textbf{Definition.} \textit{We will say that $f$ belongs to the class $\mathcal{G}$ if $f$ possesses a Gevrey-1 expansion in some right half-plane $\Re{z}>\lambda=\lambda(f)$.}

\smallskip

Let us list some properties of the class $\mathcal{G}$ required in the sequel. The proofs are either straightforward from the above definition or are given reference to.

{\bf{Property~1 (linearity).}} If $f\in\mathcal{G}$ and $g\in\mathcal{G}$, then $f+g\in\mathcal{G}$ and $af(z)\in\mathcal{G}$ for arbitrary complex $a\ne0$.

{\bf{Property~2 (shifting and dilating invariance).}} If $f\in\mathcal{G}$, then $f(Az+a)\in\mathcal{G}$ for $A>0$ and arbitrary complex  $a$.

{\bf{Property~3 (invariance under taking exponential)}} If $f(z)\in \mathcal{G}$, then $\exp\{f(z)\}\in\mathcal{G}$ (see \cite[pp. 288, 293]{Watson11} for a proof).

{\bf{Property~4.}} If $f\in\mathcal{G}$ and $\alpha_0=0$ in (\ref{eq:asimpt}), then $zf(z)\in \mathcal{G}$.

{\bf{Property~5.}} If $f\in\mathcal{G}$, then $f(z)/z\in\mathcal{G}.$

{\bf{Property~6.}} If $\varphi(z)$ is holomorphic in the neighborhood of $z=0$, then $f(z)=\varphi(a/z)\in\mathcal{G}$ for arbitrary complex $a\ne0$.

To demonstrate the last property apply Cauchy estimates to the Taylor coefficients $a_k$ of $\varphi(z)$ to get $|a_k|\leq{M}/r^k$, where $r$ is strictly less than the radius of convergence
of the Taylor series of $\varphi(z)$ at $z=0$ and $M=\max\limits_{|w|\le{r}}|\varphi(w)|$.  This leads to the following estimate of the Taylor remainder
$$
|Q_n(z)|=\left|\sum_{k=n}^{\infty}a_kz^k\right|\le\frac{M|z|^n}{r^n}\sum_{j=0}^{\infty}|z/r|^{j}=\frac{M|z|^n}{r^n(1-|z/r|)}\le\frac{M|z|^n}{r^n(1-\beta)}
$$
for $|z|\le\beta{r}$ and arbitrary $\beta\in(0,1)$.  Hence,
$$
\varphi(a/z)=\sum_{k=0}^{n-1}a_k(a/z)^k+Q_n(a/z)~\text{and}~|Q_n(a/z)|\le\frac{M|a/r|^n}{|z|^n(1-\beta)}~\text{for}~|z|\ge\frac{|a|}{\beta{r}},
$$
so that (\ref{eq:remainder}) is trivially satisfied.

We will need the next well-known lemma relating the coefficients of an asymptotic series of a function with those of its exponential. Essentially, the result contained in this lemma appeared in \cite[Appendix]{Nair}. Later, an independent derivation was given in \cite[Lemma~1]{Kalinin}. It has also been discussed recently in \cite{QSL}, where further references are given. Surprisingly, references  \cite{Nair} and  \cite{Kalinin} do not appear in \cite{QSL}.
\begin{lemma}\label{lm:exp}
Suppose $g(z)\sim\sum_{k=1}^{\infty}u_kz^{-k}$ as $z\to\infty$. Then $\exp\{g(z)\}\sim\sum_{r=0}^{\infty}v_rz^{-r}$, where the coefficients are found from
$$
v_0=1,~~~~v_r=\frac{1}{r}\sum_{k=1}^{r}ku_kv_{r-k}.
$$
Alternatively,
$$
v_r=\sum\limits_{\substack{k_1+2k_2+\cdots+rk_r=r\\k_i\ge0}}\frac{u_1^{k_1}u_2^{k_2}\cdots u_r^{k_r}}{k_1!k_2!\cdots k_r!}
=\sum\limits_{n=1}^{r}\frac{1}{n!}\sum\limits_{\substack{k_1+k_2+\cdots+k_n=r\\k_i\ge1}}\prod\limits_{i=1}^{n}u_{k_i}.
$$
\end{lemma}

\textbf{Remark.} Nair \cite[section~8]{Nair} found a determinantal expression for $v_r$ which in our notation takes the form
$$
v_r=\frac{\det(\Omega_r)}{r!},~~~\Omega_r=[\omega_{i,j}]_{i,j=1}^{r},~~\omega_{i,j}\!=\!\left\{\!\!\begin{array}{ll}u_{i-j+1}(i-j+1)(i-1)!/(j-1)!, & i\ge{j},\\
-1, &i=j-1,\\0, &i<j-1.
\end{array}\right.
$$

Various forms of the next classical theorem can be found in \cite[Theorem~4.136]{Cost}, \cite[Thereom~2.2]{DR}, \cite[paragraph~6]{GG}, \cite[section~5.7.3]{Mitschi-Sauzin}, \cite{Nevanlinna} and \cite{Sokal}.

\begin{theorem}[Watson-Nevanlinna-Sokal]\label{th:Nevanlinna}
Suppose $f(z)$ is holomorphic in $C_r=\{z:\Re(1/z)>r^{-1}\}$ and can be written as
\begin{equation}\label{eq:asim}
f(z)=\sum\limits_{k=0}^{n-1}\alpha_kz^k+R_n(z)
\end{equation}
with the error term satisfying $|R_n(z)|\leq M\tau^nn!|z|^n$, where $M$ is independent of $n$ and $z\in{C_r}$.  Then its Borel transform
\begin{equation*}
B(t)=\sum\limits_{k=0}^{\infty}\frac{\alpha_kt^k}{k!}
\end{equation*}
converges for $|t|<1/\tau$ and can be extended analytically to the domain $S_\tau=\{t: dist(t,\R^+)<1/\tau\}$ to a function satisfying
\begin{equation}\label{eq:est}
|B(t)|<K\exp(|t|/r)
\end{equation}
for some positive $K$.  Furthermore, $f$ can be recovered by the \emph{(}convergent\emph{)} integral
\begin{equation}\label{eq:int}
f(z)=\frac{1}{z}\int\limits_0^{\infty}e^{-t/z}B(t)dt.
\end{equation}
Conversely, if $B(t)$ is holomorphic in $S_{\tau'}$ \emph{(}$\tau'<\tau$\emph{)} and satisfies \emph{(\ref{eq:est})},
then the function $f(z)$ defined by the integral  \emph{(\ref{eq:int})} is holomorphic in $C_r$ and has Gevrey-1 a
symptotic approximation \emph{(\ref{eq:asim})} with uniform error bound in $C_r$,  where $\alpha_k=B^{(k)}(t)\mid_{t=0}$.
\end{theorem}

The next theorem can be found in \cite[Theorem~VIII]{Norlund14}, \cite[p.\:267]{Norlund24} and \cite[Theorem~46.2]{Wasow}.

\begin{theorem}[N{\o}rlund]\label{th:Norlund}
Suppose $\Omega(z)$ satisfies the following conditions\emph{:}
\smallskip

\emph{1)} $\Omega(z)=a/z+v(z)/[z(z+1)]$, where $v(z)$ is holomorphic and bounded in some right half-plane $\Re{z}>\kappa$.
\smallskip

\emph{2)} $\Omega(z)=\int_{0}^{\infty}e^{-tz}B(t)dt$, where $B(t)$ is holomorphic in the domain $S_\eta=\{t: dist(t,\R^+)<\eta\}$ for some $\eta>0$ and satisfies in $S_\eta$ the condition $\lim\limits_{t\to\infty}e^{-kt}B(t)=0$ for some positive $k$.

Then  $\Omega(z)$ can be expanded in the inverse factorial series
$$
\Omega(z)=\sum\limits_{s=0}^{\infty}\frac{b_s}{z(z+1)\cdots(z+s)}
$$
convergent in some right half-plane $\Re{z}>\lambda$ excluding the points $z=0,-1,-2,\ldots$.
\end{theorem}

To determine the abscissa of convergence $\lambda$ we need the following notion due to Hadamard \cite[pp.333-334]{Norlund14}:
a function $f(z)=a_0+a_1z+\cdots$  holomorphic in the unit disk $|z|<1$ has the order $\varkappa$ on the circle $|z|=1$ if
$$
\varkappa=\limsup\limits_{n\to\infty}\frac{\log|na_n|}{\log(n)}.
$$

The next theorem \cite[Theorem~III]{Norlund14} relates the order of the so-called generating function $\varphi$ of an inverse factorial series $\Omega$ to its abscissa of convergence.

\begin{theorem}[N{\o}rlund]\label{th:abscissa}
Suppose the next representation holds for sufficiently large values of $\Re{z}$\emph{:}
\begin{equation}\label{eq:seriass}
\Omega(z)=\sum\limits_{s=0}^{\infty}\frac{b_s}{(z)_{s+1}}=\int\nolimits_{0}^{1}t^{z-1}\varphi(t)dt,
\end{equation}
and assume that the order of $t^{-1}\varphi(t)$ on the circle $|1-t|=1$ is equal to $\varkappa$. If $\varkappa>1$, then the abscissa of convergence $\lambda$ of the inverse factorial series in \emph{(\ref{eq:seriass})} is equal to $\varkappa-1$, otherwise $\lambda\le\varkappa-1$.
\end{theorem}

Further, N{\o}rlund showed  in \cite[(7), page~339]{Norlund14} that a function $f(t)$ holomorphic in $|1-t|\le1$ except for a singularity at $t=0$ and representable in the form
\begin{equation}\label{eq:Norlund-order}
f(t)=\sum\limits_{i=1}^{m}t^{a_i}(\psi_{i,0}(t)+\psi_{i,1}(t)\log(t)+\cdots+\psi_{i,r}(t)\log^r(t))
\end{equation}
in the neighborhood of $t=0$  with $\Re(a_1)\le\Re(a_2)\le\cdots\le\Re(a_m)$ has the order $\kappa=-\Re(a_1)$. It is assumed that $\psi_{i,j}(t)$, $i=1,\ldots,m$, $j=0,\ldots,r$, are holomorphic around $t=0$ and such that for each $i=1,\ldots,m$ at least one of the numbers $\{\psi_{i,0}(0),\ldots,\psi_{i,r}(0)\}$ is different from zero.  Integer nonnegative $a_i$ such that $r=0$ (no logarithmic terms) must be excluded from the determination of order.

Our main tool is the following theorem.
\begin{theorem}\label{thm:seria}
Let $f(z)$ be holomorphic in some right half-plane $\Re{z}>\kappa$ and suppose that $\log(zf(z))\in\mathcal{G}$. Then for any complex  $\beta$ the function $f(z)$
can be expanded in the inverse factorial series
$$
f(z)=\sum\limits_{s=0}^{\infty}\frac{d_s}{(z+\beta)_{s+1}}
$$
convergent in some right half-plane $\Re{z}>\lambda$.
\end{theorem}

\textbf{Proof}. By Properties 3, 2 and 5, respectively, ${\zeta}f(\zeta)\in\mathcal{G}$, $(\zeta-\beta)f(\zeta-\beta)\in\mathcal{G}$ and $f(\zeta)\in\mathcal{G}$.
Then $\zeta{f(\zeta-\beta)}=(\zeta-\beta)f(\zeta-\beta)+{\beta}f(\zeta-\beta)\in\mathcal{G}$ by Properties 2 and 1. Therefore,
\begin{equation}\label{eq:asimptot2}
{\zeta}f(\zeta-\beta)=\sum\limits_{k=0}^{n-1}\frac{\beta_k}{\zeta^k}+R_n(\zeta)
\end{equation}
with the remainder $R_n(\zeta)$ bounded according to (\ref{eq:remainder}).  Now put $f_1(\zeta)=f(\zeta-\beta)$ and rewrite (\ref{eq:asimptot2}) as
\begin{equation*}\label{eq:asimpt3}
\frac{f_1(1/w)}{w}=\sum\limits_{k=0}^{n-1} {\beta_k}{w^k}+R_n(1/w)
\end{equation*}
with the error  bound of the form
$$
|R_n(1/w)|\leq {M\tau^{n}n!}{|w|^n}.
$$
According to Theorem~\ref{th:Nevanlinna}
\begin{equation*}\label{eq:int2}
\frac{f_1(1/w)}{w}=\frac{1}{w}\int\limits_0^{\infty}e^{-t/w}B(t)dt
\end{equation*}
or
\begin{equation*}\label{eq:int3}
f_1(\zeta)=\int\limits_0^\infty e^{-t\zeta}B(t)dt
\end{equation*}
and $B(t)$ satisfies condition~2 of Theorem~\ref{th:Norlund}. As $(\zeta-\beta)f_1(\zeta)\in\mathcal{G}$, we have by definition of $\mathcal{G}$:
$$
(\zeta-\beta)f_1(\zeta)=\gamma_0+\frac{\gamma_1}{\zeta}+T_2(\zeta),
\\
|T_2(\zeta)|\leq \frac{A}{|\zeta|^2}.
$$
This implies that
$$
f_1(\zeta)=\frac{\gamma_0}{\zeta}+\frac{v(\zeta)}{\zeta(\zeta+1)},
$$
where
$$
v(\zeta)=\gamma_0\beta +\gamma_1+\frac{(\gamma_0\beta +\gamma_1)(1+\beta)}{\zeta-\beta}+\frac{\zeta(\zeta+1)T_2(\zeta)}{\zeta-\beta}.
$$
The above estimate for $T_2(\zeta)$ immediately leads to conclusion that $v(\zeta)$ is bounded outside of some neighborhood of $\zeta=\beta$. Thus, the first condition of Theorem~\ref{th:Norlund} is also satisfied and
$$
f_1(\zeta)=\sum\limits_{s=0}^{\infty}\frac{d_{s}}{(\zeta)_{s+1}}.
$$
Substituting back $z=\zeta+\beta$ yields
$$
f(z)=f_1(z+\beta)=\sum\limits_{s=0}^{\infty}\frac{d_{s}}{(z+\beta)_{s+1}}.
$$
The claim regarding convergence follows by Theorem~\ref{th:Norlund}.$\hfill\square$

We will write $\Be_j(x)$  for the $j$-th Bernoulli polynomial \cite[24.2.3]{NIST}, defined by the generating function
$$
\frac{te^{xt}}{e^t-1}=\sum\limits_{n=0}^{\infty}\Be_n(x)\frac{t^n}{n!}, ~~
~~|t|<2\pi.
$$
Further, $\Be_j=\Be_j(0)$ are Bernoulli numbers  \cite[24.2.1]{NIST}.

\begin{lemma}{\label{th:Nemes}}
The function $P_a(z)=\log\Gamma(z+a)-(z+a-1/2)\log{z}+z$ belongs to the class $\mathcal{G}$ and
\begin{equation}\label{eq:Hermite-Barnes2}
P_a(z)\sim\frac{1}{2}\log(2\pi)+\sum\limits_{j=2}^{\infty}\frac{(-1)^j\Be_{j}(a)}{j(j-1)z^{j-1}}
\end{equation}
as $|z|\to\infty$ in the domain $|\arg z|<\pi-\delta$,  $0<\delta<\pi$.
\end{lemma}

\textbf{Proof.} Formula (\ref{eq:Hermite-Barnes2}) is Hermite's asymptotic expansion for $\log\Gamma(z+a)$ \cite[(1.8)]{Nemes}. We only need to proof its Gevrey-1 character.  According to \cite[(2.1.1) and (2.1.5)]{ParKam} and \cite[(10),(14)]{GG}, the Binet function
\begin{equation*}
P_0(z)-\frac{1}{2}\log(2\pi)=\log\Gamma(z)-(z-1/2)\log{z}+z-\frac{1}{2}\log(2\pi)
\end{equation*}
satisfies the relation
$$
P_0(z)-\frac{1}{2}\log(2\pi)=\sum\limits_{r=1}^{n-1}\frac{\Be_{2r}}{2r(2r-1)z^{2r-1}}+R_{2n-1}(z),
$$
in the domain $|\arg z|<\pi$ and for $z=|z|e^{i\theta}$ with $|\theta|<\pi$ the next inequality holds:
$$
|R_{2n-1}(z)|\leq\frac{|\Be_{2n}||\sec^{2n}(\theta/2)|}{2n(2n-1)|z|^{2n-1}}.
$$
In view of the asymptotic equality \cite[Corollary~1]{Dilcher}, \cite[(16)]{GG}
$$
|\Be_{2n}|=\frac{2 (2n)!}{(2\pi)^{2n}}(1+o(1)),\ n\to\infty,
$$
we conclude that $P_0(z)\in\mathcal{G}$.  It is straightforward to check that
$$
P_a(z)=P_0(z+a)-a+z\log\left(1+\frac{a}{z}\right)+\left(a-\frac{1}{2}\right)\log\left(1+\frac{a}{z}\right).
$$
Note that $P_0(z+a)\in\mathcal{G}$, $\log(1+\frac{a}{z})\in\mathcal{G}$, $z\log(1+\frac{a}{z})\in\mathcal{G}$ according to the Properties~2, 6 and 4, respectively.  It is left to apply Property~1 to conclude that $P_a(z)\in\mathcal{G}$. Uniqueness of the asymptotic expansion completes the proof. $\hfill\square$

The next lemma contains a corrected version of a formula contained in \cite{GuptaTang} which, in a different form, was already presented in \cite{Box} (see also \cite[8.5.1]{And}).
\begin{lemma}\label{lm:GammaPoincareExp}
Let $A_i,B_j>0$ satisfy $\sum_{i=1}^{p}A_i=\sum_{j=1}^{q}B_j$ and $a_i$, $b_j$, $i=1,\ldots,p$, $j=1,\ldots,q$, be arbitrary complex numbers. Then for each $0<\delta<\pi$ the next asymptotic relation holds in the domain $|\arg{z}|<\pi-\delta$,
\begin{equation*}
\frac{z^\mu}{\nu\rho^z}\frac{\prod\nolimits_{k=1}^{p}\Gamma(A_kz+a_k)}{\prod\nolimits_{j=1}^{q}\Gamma(B_jz+b_j)}\sim\sum\limits_{r=0}^\infty\frac{C_r}{z^r},
\end{equation*}
where
\begin{equation}\label{eq:nu}
\nu=(2\pi)^{(p-q)/2}\prod\nolimits_{k=1}^{p}A_k^{a_k-1/2}\prod\nolimits_{j=1}^{q}B_j^{1/2-b_j},
\end{equation}
\begin{equation}\label{eq:rho}
\rho=\prod\limits_{k=1}^{p}A_k^{A_k}\prod\limits_{j=1}^{q}B_j^{-B_j},
~~~\mu=\mu(\a,\b)=\sum\limits_{j=1}^{q}b_j-\sum\limits_{k=1}^{p}a_k+\frac{p-q}{2}.
\end{equation}
The coefficients $C_r=C_r(\A,\B;\a,\b)$ are found from the recurrence
\begin{equation}\label{eq:lr}
C_0=1,~~~~C_r=\frac{1}{r}\sum\limits_{m=1}^{r}Q_m(\A,\B;\a,\b)C_{r-m}
\end{equation}
or by other expressions contained in Lemma~\ref{lm:exp}.  Here
\begin{equation}\label{eq:qt}
Q_m(\A,\B;\a,\b)=\frac{(-1)^{m+1}}{m+1}\left[\sum\limits_{k=1}^p\frac{\Be_{m+1}(a_k)}{A_k^m}-\sum\limits_{j=1}^q\frac{\Be_{m+1}(b_j)}{B_j^m}
\right],
\end{equation}
and $\Be_{m}(\cdot)$ denotes the $m$-th Bernoulli polynomial. Furthermore,
$$
\log\left(\frac{z^\mu}{\nu\rho^z}\frac{\prod\nolimits_{k=1}^{p}\Gamma(A_kz+a_k)}{\prod\nolimits_{j=1}^{q}\Gamma(B_jz+b_j)}\right)\in\mathcal{G}.
$$
\end{lemma}
\textbf{Proof.} From Hermite's asymptotic formula (\ref{eq:Hermite-Barnes2}) we obtain by straightforward calculations:
\begin{multline*}
\log\Gamma(Az+a)=(Az+a-1/2)\log{Az}-Az+\frac{1}{2}\log(2\pi)+\sum\limits_{j=2}^{m}\frac{(-1)^j\Be_{j}(a)}{j(j-1)(Az)^{j-1}}+R_{m}(z)
\\
=Az\log{z}+(A\log{A}-A)z+(a-1/2)\log{z}+(a-1/2)\log{A}+\frac{1}{2}\log(2\pi)
\\
+\sum\limits_{j=2}^{m}\frac{(-1)^j\Be_{j}(a)}{j(j-1)A^{j-1}z^{j-1}}+R_{m}(z),
\end{multline*}
where  $R_m(z)=\O(z^{-m})$ in the sector $|\arg(z)|<\pi-\delta$ for any $\delta>0$. For $A>0$ Lemma~\ref{th:Nemes} implies that
$\log\Gamma(Az+a)-(Az+a-1/2)\log{Az}+Az\in\mathcal{G}$.  Then employing the condition $\sum{A_k}=\sum{B_j}$ we obtain
\begin{equation*}
\log\left\{\frac{\prod\nolimits_{k=1}^{p}\Gamma(A_kz+a_k)}{\prod\nolimits_{j=1}^{q}\Gamma(B_jz+b_j)}\right\} \sim z\log{\rho}-\mu\log{z}+\log{\nu}
+\sum\limits_{t=1}^{\infty}\frac{Q_t(\A,\B;\a,\b)}{tz^{t}},
\end{equation*}
where $Q_t(\A,\B;\a,\b)$ is given in (\ref{eq:qt}). According to Property~1
\begin{equation*}
\log\left\{\frac{z^\mu}{\nu\rho^z}\frac{\prod\nolimits_{k=1}^{p}\Gamma(A_kz+a_k)}{\prod\nolimits_{j=1}^{q}\Gamma(B_jz+b_j)}\right\}\in\mathcal{G}.
\end{equation*}
Formula (\ref{eq:lr}) and other methods to compute $C_r$ follow from Lemma~\ref{lm:exp}. $\hfill\square$

\paragraph{3. Main results.}
To formulate our mains theorem we will need the non-central Stirling numbers of the first kind \cite[8.5]{Charalambides} defined by
\begin{equation}\label{eq:noncentralStirling}
s_{\sigma}(n,r)=\sum\limits_{k=r}^{n}(-1)^{k+r}\binom{n}{k}(\sigma)_{n-k}s(k,r)=\sum\limits_{k=0}^{n}(-1)^{k+r}\binom{n}{k}(\sigma)_{n-k}s(k,r),
\end{equation}
where $s(n,j)$ denotes the ordinary Stirling number of the first kind  generated by
$$
x(x-1)\cdots(x-n+1)=\sum_{j=0}^{n}s(n,j)x^j.
$$
Their ''horizontal'' exponential generating function is given by \cite[(A.2)]{Weniger2010}:
\begin{multline*}
\sum\limits_{l=0}^{n}x^{l}s_{\sigma}(n,l)=\sum\limits_{l=0}^{n}x^{l}\sum\limits_{k=0}^{n}(-1)^{k+l}\binom{n}{k}(\sigma)_{n-k}s(k,l)
\\
=\sum\limits_{k=0}^{n}(-1)^k\binom{n}{k}(\sigma)_{n-k}\sum\limits_{l=0}^{k}s(k,l)(-x)^{l}
=\sum\limits_{k=0}^{n}\binom{n}{k}(\sigma)_{n-k}(x)_{k}=(\sigma+x)_{n},
\end{multline*}
where the expansion  $(x)_{k}=(-1)^k\sum_{l=0}^{k}s(k,l)(-x)^{l}$ has been applied.
Sometimes it is more convenient to use the  ''vertical''  generating function
\begin{multline*}
\sum\limits_{n=0}^{\infty}\frac{s_{\sigma}(n,l)}{n!}x^{n}=\sum\limits_{n=l}^{\infty}\frac{s_{\sigma}(n,l)}{n!}x^{n}
=\sum\limits_{n=l}^{\infty}x^{n}\sum\limits_{k=l}^{n}(-1)^{k+l}\frac{(\sigma)_{n-k}s(k,l)}{k!(n-k)!}
\\
=(-1)^{l}\sum\limits_{k=l}^{\infty}(-1)^{k}s(k,l)\frac{x^{k}}{k!}\sum\limits_{n=k}^{\infty}\frac{(\sigma)_{n-k}}{(n-k)!}x^{n-k}
\\
=(-1)^{l}\sum\limits_{k=l}^{\infty}s(k,l)\frac{(-x)^{k}}{k!}\sum\limits_{m=0}^{\infty}\frac{(\sigma)_{m}}{m!}x^{m}
=\frac{1}{l!}\frac{1}{(1-x)^{\sigma}}\left(\log\frac{1}{1-x}\right)^{\!l}.
\end{multline*}

The non-central Stirling numbers were studied by Carlitz in \cite{Carlitz1,Carlitz2} using the symbol $R_1(n,l,\sigma)=s_{\sigma}(n,l)$  and some years later also by
Koutras \cite{Koutras}. Broder \cite{Broder} considered them for integer $\sigma$ from the combinatorial viewpoint. Various formulas for these numbers are given in  \cite[8.5]{Charalambides}. Among other things, Carlitz found the double generating function \cite[(5.4)]{Carlitz1}
$$
\sum\limits_{l,n=0}^{\infty}s_{\sigma}(n,l)y^l\frac{x^n}{n!}=(1-x)^{-\sigma-y}.
$$
Using a formula from \cite[section~6.43,p.134]{Milne-Thompson} this generating function leads to the connection formula  \cite[(7.6)]{Carlitz2}
$$
s_{\sigma}(n,l)=\binom{-l-1}{n-l}\Be^{(n+1)}_{n-l}(1-\sigma)=\frac{(-1)^{n-l}(l+1)_{n-l}}{(n-l)!}\Be^{(n+1)}_{n-l}(1-\sigma),
$$
where $\Be^{\gamma}_{k}(x)$ is the $k$-th N{\o}rlund-Bernoulli polynomial (also known as the generalized Bernoulli polynomials) defined by the generating function \cite[(1)]{Norlund61}
\begin{equation}\label{eq:NorlundBernoulli}
\frac{t^{\gamma}e^{xt}}{(e^t-1)^{\gamma}}=\sum\limits_{k=0}^{\infty}\Be^{(\gamma)}_{k}(x)\frac{t^k}{k!}.
\end{equation}

The main result of this note is the following theorem.
\begin{theorem}\label{th:seria2}
Let $A_k,B_j>0$ satisfy $\sum_{k=1}^{p}A_k=\sum_{j=1}^{q}B_j$ and $\sigma$, $a_k$, $b_j$, $k=1,\ldots,p$, $j=1,\ldots,q$, be  complex numbers such that $\mu=1$, where $\mu$ is defined in \emph{(\ref{eq:rho})}. Suppose further that $\nu$ and $\rho$ are  given by \emph{(\ref{eq:nu})} and  \emph{(\ref{eq:rho})}, respectively.  Then the inverse factorial series expansion
\begin{equation}\label{eq:GammaRatioExt4}
W(z)=\rho^{-z}\frac{\prod\nolimits_{k=1}^{p}\Gamma(A_kz+a_k)}{\prod\nolimits_{j=1}^{q}\Gamma(B_jz+b_j)}
=\sum\limits_{n=0}^{\infty}\frac{\nu}{(z+\sigma)_{n+1}}\sum\limits_{r=0}^{n}C_{r}s_{\sigma}(n,r)
\end{equation}
converges for $\Re{z}>\alpha$ except at the points $z=-\sigma-l$,  $l\in\N_0$. Here $\alpha$ denotes the real part of the rightmost pole of $W(z)$ \emph{(}obviously, $\alpha=-\min\{\Re(a_k/A_k): k=1,\ldots,p\}$ if no cancelations of the numerator and denominator poles take place\emph{)}.  Moreover, if $\Re(\sigma+\alpha)>0$, then $\alpha$ is equal to the abscissa of convergence of the series \emph{(\ref{eq:GammaRatioExt4})}.   Here $C_{r}=C_{r}(\A,\B;\a,\b)$ is defined in \emph{(\ref{eq:lr})}, \emph{(\ref{eq:qt})} and $s_{\sigma}(n,r)$ is the non-central Stirling number of the first kind.
\end{theorem}
\textbf{Remark.} The convergence domain in the above theorem can be described as follows: if the abscissa of the rightmost pole of $W(z)$ is greater than the abscissa of the rightmost pole of the series on the right hand side, then the convergence abscissa is equal precisely to the abscissa of the rightmost pole of $W(z)$; if the abscissa of the rightmost pole of $W(z)$ does not exceed the abscissa of the rightmost pole of the series on the right hand side, then the convergence abscissa does not exceed the abscissa of the rightmost pole of $W(z)$; if the rightmost poles on both sides are simple and coincide (i.e. $\sigma=a_{k^*}/A_{k^*}$, where $k^*=\argmin\{\Re(a_k/A_k): k=1,\ldots,p\}$), then these poles should be ignored when calculating the convergence abscissa.

\textbf{Proof.} As $\mu=1$ according to Lemma~\ref{lm:GammaPoincareExp} $\log(zW(z))\in\mathcal{G}$. Then, by Theorem~\ref{thm:seria} we conclude that $W(z)$ can be expanded in a series of inverse factorials $[(z+\sigma)_{n+1}]^{-1}$.  As this series is also asymptotic for $W(z)$ as $z\to\infty$, its coefficients can be obtained by rearranging the Poincar\'{e} asymptotic expansion of $W(z)$.  According to Lemma~\ref{lm:GammaPoincareExp} the latter is given by
$$
W(z)\sim\nu\sum\limits_{r=0}^{\infty}\frac{C_{r}(\A,\B;\a,\b)}{z^{r+1}},
$$
where the coefficients  $C_{r}(\A,\B;\a,\b)$ are defined in (\ref{eq:lr}), (\ref{eq:qt}). Following \cite{Watson12} and  \cite[(4.1)]{Weniger2010} this asymptotic series can be re-expanded using the asymptotic sequence $\{1/(z)_{n+1}\}_{n=0}^{\infty}$  by applying \cite[\S30(6)]{Nielsen}
$$
\frac{1}{z^{r+1}}=\sum\limits_{k=0}^{\infty}\frac{(-1)^ks(r+k,r)}{(z)_{r+k+1}},
$$
where $s(n,j)$ is the Stirling number of the first kind. Substituting and changing the order of summations we get
$$
W(z)\sim\nu\sum\limits_{n=0}^{\infty}\frac{(-1)^n}{(z)_{n+1}}\sum\limits_{r=0}^{n}(-1)^rs(n,r)C_{r}(\A,\B;\a,\b).
$$
Next, instead of inverse factorials $1/(z)_{m}$ we can utilize the asymptotically equivalent sequence
$\{1/(z+\sigma)_{m}\}_{m\ge0}$ by employing the connection formula (see \cite[(10)]{Norlund14} or \cite[(138.15)]{Norlund24})
\begin{equation}\label{eq:shiftconnection}
U(z)=\sum\limits_{k=0}^{\infty}\frac{u_{k+1}k!}{(z)_{k+1}}=\sum\limits_{k=0}^{\infty}\frac{v_{k+1}k!}{(z+\sigma)_{k+1}},
\end{equation}
with the  coefficients related by
\begin{equation}\label{eq:ak-bktransform}
v_{k+1}=\sum\limits_{j=0}^{k}\binom{\sigma+j-1}{j}u_{k+1-j}=\sum\limits_{j=0}^{k}\frac{(\sigma)_j}{j!}u_{k+1-j}.
\end{equation}
The inverse formula reads \cite[2.1(1)]{Riordan}
$$
u_{k+1}=\sum\limits_{j=0}^{k}(-1)^{j}\binom{\sigma}{j}v_{k+1-j}.
$$
Combining these facts with (\ref{eq:noncentralStirling}), we arrive at ($C_j=C_j(\A,\B;\a,\b)$)
\begin{multline*}
\frac{W(z)}{\nu}\!=\!\!\sum\limits_{n=0}^{\infty}\!\frac{(-1)^n}{(z)_{n+1}}\sum\limits_{j=0}^{n}(-1)^js(n,j)C_j
\!=\!\!\sum\limits_{n=0}^{\infty}\!\frac{n!}{(z+\sigma)_{n+1}}\sum\limits_{j=0}^{n}\!\frac{(\sigma)_j(-1)^{n-j}}{(n-j)!j!}\sum\limits_{r=0}^{n-j}(-1)^rs(n-j,r)C_r
\\
=\sum\limits_{n=0}^{\infty}\frac{1}{(z+\sigma)_{n+1}}\sum\limits_{r=0}^{n}(-1)^rC_r\sum\limits_{j=0}^{n-r}(-1)^{n-j}\binom{n}{j}(\sigma)_js(n-j,r)
\\
=\sum\limits_{n=0}^{\infty}\frac{1}{(z+\sigma)_{n+1}}\sum\limits_{r=0}^{n}(-1)^rC_r
\sum\limits_{k=r}^{n}(-1)^{k}\binom{n}{k}(\sigma)_{n-k}s(k,r)=\sum\limits_{n=0}^{\infty}\frac{1}{(z+\sigma)_{n+1}}\sum\limits_{r=0}^{n}C_{r}s_{\sigma}(n,r)
\end{multline*}
which is precisely (\ref{eq:GammaRatioExt4}).

To compute the convergence abscissa denote $x=z+\sigma$ and consider the function
$$
\Omega(x)=\rho^{-\sigma}W(z)=\rho^{-x}\frac{\prod\nolimits_{k=1}^{p}\Gamma(A_kx+a_k-A_k\sigma)}{\prod\nolimits_{j=1}^{q}\Gamma(B_jx+b_j-B_j\sigma)}.
$$
Denote $\a'=\a-\sigma\A$, $\b'=\b-\sigma\B$.  Then, clearly, $\mu(\a',\b')=\mu(\a,\b)=1$ (by the hypotheses of the theorem; $\mu(\a,\b)$ is defined in \ref{eq:rho}).
According to \cite[Theorem~6]{KPCMFT2016}
$$
\Omega(x)=\int_0^1 t^{x-1}H_{q,p}^{p,0}\left(\rho{t}\left|\!\!\begin{array}{c}(\B,\b')\\(\A,\a')\end{array}\right.\!\right)dt,
$$
where $H_{q,p}^{p,0}(\rho{t})$ is a particular case of Fox's $H$-function defined by
\begin{equation*}
H_{q,p}^{p,0}\left(\rho{t}\left|\!\!\begin{array}{c}(\B,\b')\\(\A,\a')\end{array}\right.\!\right)=\frac{1}{2\pi{i}}\int\limits_{c-i\infty}^{c+i\infty}
\frac{\prod\nolimits_{k=1}^{p}\Gamma(A_ks+a_k-A_k\sigma)}{\prod\nolimits_{j=1}^{q}\Gamma(B_j s+b_j-B_j\sigma)}(\rho{t})^{-s}ds
\end{equation*}
(further details about the definition of Fox's $H$ function can be found in \cite[Chapter~1]{KilSaig}; see also \cite{KPCMFT2016,KPCMFT2017}). By \cite[Theorem~1.5]{KilSaig} the following representation is true
\begin{equation*}
H_{q,p}^{p,0}\left(\rho{t}\left|\!\!\begin{array}{c}(\B,\b')\\(\A,\a')\end{array}\right.\!\right)=
\sum\limits_{k=1}^{p}\sum\limits_{l=0}^{\infty}\sum_{i=0}^{N_{kl}-1}H_{kli} (\rho{t})^{(a_k-A_k\sigma+l)/A_k}\log^i(\rho{t}),
\end{equation*}
where the two outer summations run over all poles of the function $\Omega(x)$ and $N_{kl}$ denotes the order of the pole at the point $a_{kl}=(-a_k+A_k\sigma-l)/A_k$.  The explicit form of the constants $H_{kli}$ is immaterial here. We emphasize, however, that if the pole of the numerator at $x=a_{kl}$ cancels out with a pole of the denominator, then we put $N_{kl}=0$ and the corresponding term is omitted from the above summation.
Comparing this representation with (\ref{eq:Norlund-order}) we conclude by N{\o}rlund's argument \cite[p.339]{Norlund14} explained below (\ref{eq:Norlund-order}) that
the order $\varkappa$ of $H_{q,p}^{p,0}(\rho{t})/t$ on the circle $|1-t|=1$ is given by
$$
\varkappa=-\min_{k}\left\{\Re(a_k/A_k-\sigma-1)\right\}=\max_{k}\left\{\Re(-a_k/A_k)\right\}+\Re(\sigma)+1,
$$
where the minimum is taken over the indices $k\in\{1,\ldots,p\}$, such that $\Omega(x)$ has a pole at $x_k=-a_k/A_k+\sigma$ and $-x_k$ is not a nonnegative integer.  Then by Theorem~\ref{th:abscissa} the abscissa of convergence $\lambda_{\Omega}$ of $\Omega$ satisfies $\lambda_{\Omega}\le\varkappa-1=\max_{k}\left\{\Re(-a_k/A_k)\right\}+\Re(\sigma)$ if $\varkappa-1\leq0$ or  $-\min_{k}\left\{\Re(a_k/A_k)\right\}\le-\Re(\sigma)$. The last condition can be interpreted as follows: the abscissa of the rightmost pole of $W(z)$ does not exceed the abscissa of the rightmost pole of the series on the right hand side of (\ref{eq:GammaRatioExt4}).   So that by the reverse change of variable $z=x-\sigma$ the convergence abscissa $\lambda$ of the original series (\ref{eq:GammaRatioExt4}) satisfies $\lambda\le\max_{k}\left\{\Re(-a_k/A_k)\right\}$.

If, on the contrary, $\varkappa-1>0$ or $-\min_{k}\left\{\Re(a_k/A_k)\right\}>-\Re(\sigma)$, then $\lambda_{\Omega}=\varkappa-1=\max_{k}\left\{\Re(-a_k/A_k)\right\}+\Re(\sigma)$ implying that the convergence abscissa $\lambda$ of the original series (\ref{eq:GammaRatioExt4}) is $\lambda=\max_{k}\left\{\Re(-a_k/A_k)\right\}$ (the abscissa of the rightmost pole of $W(z)$).  Note, that for the rightmost pole $x_{k^*}=-a_{k^*}/A_{k^*}+\sigma$ the situation $-x_{k^*}\in\N_0$ contradicts $-\min_{k}\left\{\Re(a_k/A_k)\right\}>-\Re(\sigma)$, so that under this condition we necessarily have the equality $\lambda=\max_{k}\left\{\Re(-a_k/A_k)\right\}$.  $\hfill\square$

\medskip

It is easy to modify the above theorem to get rid of the restriction $\mu=1$.
\begin{corollary}\label{cr:gammaratio6}
Let $A_k,B_j>0$ satisfy $\sum_{k=1}^{p}A_k=\sum_{j=1}^{q}B_j$ and let $\theta$, $a_k,b_j$, $k=1,\ldots,p$,  $j=1,\ldots,q$, be arbitrary  complex numbers.  Suppose further that $\nu$, $\rho$ and $\mu$ are defined in \emph{(\ref{eq:nu}) and (\ref{eq:rho})}, respectively.  Then
\begin{equation}\label{eq:GeneralW}
W(z)=\rho^{-z}\frac{\prod\nolimits_{k=1}^{p}\Gamma(A_kz+a_k)}{\prod\nolimits_{j=1}^{q}\Gamma(B_jz+b_j)}
=\sum\limits_{n=0}^{\infty}\frac{h_n\Gamma(z+\theta+1)}{\Gamma(z+\theta+\mu+n+1)},
\end{equation}
where the series converges in the half-plane $\Re{z}>\beta$ with $\beta$ equal to the real part of the rightmost pole of the function $W(z)\Gamma(z+\theta+\mu)/\Gamma(z+\theta+1)$.  The coefficients are computed by
$$
h_n={\nu}\sum\limits_{r=0}^{n}C_r(\A',\B';\a',\b')s_{\theta+\mu}(n,r)=\nu\Gamma(n+\mu)\sum\limits_{r+k=n}\frac{(-1)^kC_r(\A,\B;\a,\b)}{k!\Gamma(r+\mu)}\Be^{(n+\mu)}_{k}(-\theta),
$$
where $\A'=(\A,1)$, $\B'=(\B,1)$, $\a'=(\a,\theta+\mu)$, $\b'=(\b,\theta+1)$ and $C_r(\A',\B';\a',\b')$ are defined in \emph{(\ref{eq:lr})}.  Moreover, if $\Re(\beta+\theta+\mu)>0$, then $\beta$ is equal to the convergence abscissa of the series \emph{(\ref{eq:GeneralW})}.
\end{corollary}
\textbf{Proof.} If $\mu\neq1$ we can take $\sigma=\theta+\mu$ in Theorem~\ref{th:seria2} to get
\begin{equation*}
W_1(z)=\rho^{-z}\frac{\prod\nolimits_{k=1}^{p}\Gamma(A_kz+a_k)}{\prod\nolimits_{j=1}^{q}\Gamma(B_jz+b_j)}
\frac{\Gamma(z+\theta+\mu)}{\Gamma(z+\theta+1)}=\sum\limits_{n=0}^{\infty}\frac{h_n}{(z+\theta+\mu)_{n+1}},
\end{equation*}
where $h_n={\nu}\sum_{r=0}^{n}C_r(\A',\B';\a',\b')s_{\theta+\mu}(n,r)$ and the numbers $C_r(\A',\B';\a',\b')$ are defined in (\ref{eq:lr}).  Indeed,  $\mu(\a',\b')=1$ by construction, so that Theorem~\ref{th:seria2} is applicable for $W_1(z)$. Multiplying both sides of (\ref{eq:GammaRatioExt4}) by $\Gamma(z+\theta+1)/\Gamma(z+\theta+\mu)$ we arrive at (\ref{eq:GeneralW}) together with the first formula for $h_n$ and all claims regarding convergence.  It remains to verify the second expression for $h_n$. To this end employ the expansion \cite[(43)]{Norlund61}
$$
\frac{1}{z^{\beta}}=\sum\limits_{k=0}^{\infty}\frac{(-1)^k\Be^{(\beta+k)}_{k}(-\theta)(\beta)_{k}\Gamma(z+\theta+1)}{k!\Gamma(z+\theta+\beta+k+1)},
$$
where the series is known to converge for $\Re{z}>0$. Applying Lemma~\ref{lm:GammaPoincareExp} and making the necessary rearrangements we get
\begin{multline*}
\rho^{-z}\frac{\prod\nolimits_{k=1}^{p}\Gamma(A_kz+a_k)}{\prod\nolimits_{j=1}^{q}\Gamma(B_jz+b_j)}
\sim\nu\sum\limits_{r=0}^\infty\frac{C_r}{z^{r+\mu}}=\nu\sum\limits_{r=0}^{\infty}C_r
\sum\limits_{k=0}^{\infty}\frac{(-1)^k\Be^{(r+\mu+k)}_{k}(-\theta)(r+\mu)_{k}\Gamma(z+\theta+1)}{k!\Gamma(z+\theta+r+\mu+k+1)}
\\
=\nu\sum\limits_{n=0}^{\infty}\frac{\Gamma(z+\theta+1)\Gamma(n+\mu)}{\Gamma(z+\theta+\mu+n+1)}\sum\limits_{r+k=n}
\frac{(-1)^kC_r}{k!\Gamma(r+\mu)}\Be^{(n+\mu)}_{k}(-\theta),
\end{multline*}
where this time $C_r=C_r(\A,\B;\a,\b)$.  As the inverse factorial series of a given function is unique whether it is convergent or asymptotic we conclude that
$$
h_n=\nu\Gamma(n+\mu)\sum\limits_{r+k=n}\frac{(-1)^kC_r(\A,\B;\a,\b)}{k!\Gamma(r+\mu)}\Be^{(n+\mu)}_{k}(-\theta).~~~~~\square
$$

To conclude this section we remark that the expansions presented in this note can probably be generalized  to the case when $\Delta=\sum_{j=1}^{q}B_j-\sum_{k=1}^{p}A_k\ne0$.  For $\Delta>0$ this expansion will be in terms of reciprocals of the gamma functions $\Gamma(\Delta{z}+\sigma)$, while for $\Delta<0$ the expansion will be in terms of $\Gamma(-\Delta{z}+\sigma)$ (not reciprocals!).  To derive such expansions one may apply the technique developed by Riney in \cite{Riney2}, where the general $\Delta>0$ case is deduced from the expansion similar to (\ref{eq:GeneralW}) for $\Delta=0$.  Riney's results, however, are only asymptotic and pertain to the ''unweighted'' case $A_k=B_j=1$.

\paragraph{4. An identity for the N{\o}rlund-Bernoulli polynomials.}
In this section we present an identity for the N{\o}rlund-Bernoulli defined by the generating function (\ref{eq:NorlundBernoulli}). Although its novelty is dubious, we could not locate it in the existing literature. The identity is obtained by
comparing different expansion for the simplest particular case of the function (\ref{eq:W}) and is given in the next proposition.
\begin{theorem}\label{th:Bernulli}
The N{\o}rlund-Bernoulli polynomials satisfy
$$
\frac{(t-x+1)_{m}}{m!}\Be^{(t-x+m+1)}_{m}(1-x)
=\sum\limits_{j=0}^{m}\Be^{(t-x+1)}_{j}(t)\frac{(x-t)_j}{j!}\sum\limits_{k=0}^{m-j}(-1)^k\binom{m}{k}s(m-k,j)(t)_k
$$
and
$$
\frac{1}{m!}\sum\limits_{j=0}^{m}s(m,j)\Be^{(t+x)}_{j}(t)\frac{(1-t-x)_j}{j!}=\sum\limits_{j=0}^{m}\binom{t}{m-j}\Be^{(t+x+j)}_{j}(x)
\frac{(t+x)_j}{(j!)^2},
$$
where $s(m,j)$ is the Stirling number of the first kind.
\end{theorem}

\textbf{Proof.} We will use an asymptotic expansion in inverse powers of $z$ for the ratio $\Gamma(z+t)/\Gamma(z+x)$ found by Tricomi and Erd\'{e}lyi in \cite{TricErdelyi}.  This expansion is given on page 141 of \cite{TricErdelyi}, but we will rewrite it in a slightly different form using the definition of the coefficients \cite[(19)]{TricErdelyi} and the identity $\Gamma(1-a)/\Gamma(1-a-n)=(-1)^n(a)_n$:
\begin{equation}\label{eq:TricErd}
\frac{\Gamma(z+t)}{\Gamma(z+x)z^{t-x+1}}\sim\sum\limits_{n=0}^{\infty}\frac{(-1)^n\Be^{(t-x+1)}_{n}(t)(x-t)_n}{n!z^{n+1}}~~\text{as}~~z\to\infty,
\end{equation}
where $-x,-t,x-t\notin\N$ and $z$ is in $\C$ cut along the ray connecting $-\alpha$ and $-\alpha-\infty$. Following \cite{Weniger2010}, we can substitute the expansion \cite[(6) on p. 78]{Nielsen}
$$
\frac{1}{z^{n+1}}=\sum\limits_{j=0}^{\infty}\frac{(-1)^js(n+j,n)}{(z)_{n+j+1}}
$$
in (\ref{eq:TricErd}) to get
\begin{multline}\label{eq:invfac1}
\frac{\Gamma(z+t)}{\Gamma(z+x)z^{t-x+1}}\sim\sum\limits_{n=0}^{\infty}
\frac{(-1)^n\Be^{(t-x+1)}_{n}(t)(x-t)_n}{n!}\sum\limits_{j=0}^{\infty}\frac{(-1)^js(n+j,n)}{(z)_{n+j+1}}
\\
=\sum\limits_{m=0}^{\infty}\frac{(-1)^m}{(z)_{m+1}}\sum\limits_{n=0}^{m}\frac{s(m,n)\Be^{(t-x+1)}_{n}(t)(x-t)_n}{n!}.
\end{multline}
On the other hand N{\o}rlund derived the expansion \cite[(43)]{Norlund61}
\begin{equation}\label{eq:NorlundGammaRatio}
\frac{\Gamma(z+t)}{\Gamma(z+x)z^{t-x+1}}=\sum\limits_{m=0}^{\infty}\frac{(-1)^m\Be^{(t-x+m+1)}_{m}(1-x)(t-x+1)_m}{m!(z+t)_{m+1}}
\end{equation}
convergent in some right half-plane.  Equating the coefficients in (\ref{eq:invfac1}) and (\ref{eq:NorlundGammaRatio}) we
get an  identity of the form (\ref{eq:shiftconnection}) with
$$
v_{m+1}=\frac{(-1)^m\Be^{(t-x+m+1)}_{m}(1-x)(t-x+1)_m}{(m!)^2},
~~u_{m+1}=\frac{(-1)^m}{m!}\sum\limits_{n=0}^{m}\frac{s(m,n)\Be^{(t-x+1)}_{n}(t)(x-t)_n}{n!}.
$$
After some rearrangement  (\ref{eq:ak-bktransform}) takes the form
\begin{multline}
\frac{(t-x+1)_{m}}{m!}\Be^{(t-x+m+1)}_{m}(1-x)
=\sum\limits_{k=0}^{m}(-1)^k\binom{m}{k}(t)_k\sum\limits_{j=0}^{m-k}s(m-k,j)\Be^{(t-x+1)}_{j}(t)\frac{(x-t)_j}{j!}
\\
=\sum\limits_{j=0}^{m}\Be^{(t-x+1)}_{j}(t)\frac{(x-t)_j}{j!}\sum\limits_{k=0}^{m-j}(-1)^k\binom{m}{k}s(m-k,j)(t)_k.
\end{multline}
Applying the inversion formula below (\ref{eq:ak-bktransform}) and changing $1-x\mapsto{x}$  we finally obtain
$$
\frac{1}{m!}\sum\limits_{j=0}^{m}s(m,j)\Be^{(t+x)}_{j}(t)\frac{(1-t-x)_j}{j!}=\sum\limits_{j=0}^{m}\binom{t}{m-j}\Be^{(t+x+j)}_{j}(x)
\frac{(t+x)_j}{(j!)^2}.~~~\square
$$

\bigskip

\textbf{Acknowledgements}.  We thank Gerg\H{o} Nemes for contributing a lemma to the original version of this paper and useful discussions.  This research was supported by the Russian Science Foundation under project 14-11-00022.

\end{document}